\documentclass[11pt,a4paper]{article}
\usepackage[english]{babel}
\usepackage{amsmath}
\usepackage{amssymb}
\usepackage{graphicx}

\newcommand{\ds}{\displaystyle}
\newtheorem{de}{Definition}[section]
\newtheorem{theo}{Theorem}[section]
\newtheorem{cor}[theo]{Corollary}
\newtheorem{prop}[theo]{Proposition}
\newtheorem{lemm}[theo]{Lemma}
\newtheorem{rem}[theo]{Remark}
\newcommand{\nm}{\noalign{\smallskip}}

\newcommand{\R}{\mathbb{R}}
\newcommand{\Ra }{{\cal R}}

\newcommand{\M}{{\cal M}}
\newcommand{\A}{{\cal A}}

\newcommand{\C}{{\mathcal C}}

\newcommand{\ii}{\infty}

\newcommand{\qed}{\hfill $\Box$ \medskip}

\newcommand{\1}{\textbf{1}}
\newcommand{\Om}{\Omega}

\newcommand{\s}{\sigma}
\newcommand{\al}{\alpha}

\newcommand{\be}{\beta}
\newcommand{\bs}{\backslash}

\newcommand{\e}{\varepsilon}
\newcommand{\dr}{\partial}
\newcommand{\del}{\delta}
\newcommand{\g}{\nabla}
\newcommand{\ti}{\widetilde}

\newcommand{\ph}{\varphi}
\newcommand{\La}{\triangle}
\newcommand{\Proof}{\medskip \noindent {\sl Proof}. \ }
\newcommand{\h}{{H^1_0(\Om)}}
\newcommand{\hh}{{H^{-1}(\Om)}}

\newcommand{\fg}{\longrightarrow}
\newcommand{\cqfd}{\hfill $\square$}
\newcommand{\n}[2]{\|{#1}\|_{#2}}

\newcommand{\Div}{\nabla \cdot}

\newcommand{\Hd}{H_0^1(D)}

\numberwithin{equation}{section}

\begin{document}

\title{Reconstruction of a piecewise smooth absorption coefficient by an acousto-optic
process\thanks{\footnotesize This work was
 supported by the ERC Advanced Grant
Project MULTIMOD--267184.}}
\author{Habib Ammari\thanks{\footnotesize Department of Mathematics and Applications,
Ecole Normale Sup\'erieure, 45 rue d'Ulm - F 75230 Paris cedex 05,
France (habib.ammari@ens.fr, lnguyen@dma.ens.fr,
laurent.seppecher@ens.fr).} \and Josselin
Garnier\thanks{\footnotesize Laboratoire de Probabilit\'es et
Mod\`eles Al\'eatoires \& Laboratoire Jacques-Louis Lions,
Universit\'e Paris VII, 75205 Paris Cedex 13, France
(garnier@math.jussieu.fr).}  \and Loc Hoang Nguyen\footnotemark[2]
\and Laurent Seppecher\footnotemark[2]}

\maketitle
\begin{abstract}
The aim of this paper is to tackle the nonlinear optical
reconstruction problem. Given a set of acousto-optic measurements,
we develop a mathematical framework for the reconstruction problem
in the case where the optical absorption distribution is supposed
to be a perturbation of a piecewise constant function. Analyzing
the acousto-optic measurements, we establish a new equation in the
sense of distributions for the optical absorption coefficient. For
doing so, we introduce a weak Helmholtz decomposition and
interpret in a weak sense the cross-correlation measurements using
the spherical Radon transform. We next show how to find an initial
guess for the unknown coefficient and finally construct the true
coefficient by providing a Landweber type iteration and proving
that the resulting sequence converges to the solution of the
system constituted by the optical diffusion equation and the new
equation mentioned above. Our results in this paper generalize the
acousto-optic process proposed in
\cite{AmmariLaurentetal:preprint2012} for piecewise smooth optical
absorption distributions.
\end{abstract}

\bigskip

\noindent {\footnotesize Mathematics Subject Classification
(MSC2000): 35R30, 35B30.}


\section{Introduction} \label{sec: intro}
Let $\Om$ be a bounded $\C^1$-domain of $\R^d$, where $d \in
\{2,3\}$. We denote by $\nu$ the outward normal to $\partial
\Omega$, the boundary of $\Omega$. We need the following
functional spaces. For $m$ a non-negative integer, we define the
space $H^{m}(\Omega)$ as the family of all functions in
$L^2(\Omega)$, whose weak derivatives of orders up to $m$ also
belong to $L^2(\Omega)$. For $m\geq 1$, the space
$H^{m-1/2}(\partial \Omega)$ denotes the set of the traces on
$\partial \Omega$ of all functions in $H^m(\Omega)$. We let
$H_0^{m}(\Omega)$ be the closure of $\mathcal{C}^\infty_c(\Omega)$
in $H^{m}(\Omega)$, where $\mathcal{C}^\infty_c(\Omega)$ is the
set of all infinitely differentiable functions with compact
supports in $\Omega$. We denote by $H^{-m}(\Omega)$ the dual of
$H_0^{m}(\Omega)$. Finally, for $p\geq 1$, we introduce
$W^{m,p}(\Omega)$ as the space of functions whose weak derivatives
of orders up to $m$ are functions in $L^p(\Omega)$ and
$W_0^{m,p}(\Omega)$ to be the closure of
$\mathcal{C}^\infty_c(\Omega)$ in $W^{m,p}(\Omega)$. Here,
$L^p(\Omega)$ is defined in the usual way. Note that
$W^{m,2}(\Omega) = H^m(\Omega)$ and $W_0^{m,2}(\Omega) =
H^m_0(\Omega)$.

Suppose that $\Omega$ represents an optical medium and let $a_*:\Om
\fg \R^+$ be the optical absorption coefficient of $\Omega$. When
the medium $\Omega$ is illuminated with infrared light spots, the
optical energy density $\Phi_* \in H^2(\Omega)$ inside $\Omega$
satisfies the diffusion equation
\begin{equation}
\left\{\begin{aligned}
    -\La\Phi_* +a_*\Phi_* &=0\ \text{ in } \Om,\\
    \nm
    l \dr_\nu\Phi_* + \Phi_* &= g\ \text{ on } \partial\Om ,
\end{aligned} \right.
\label{Optical eqn}
\end{equation}
where $l \geq 0$ is the extrapolation length, computed from the
radiative transport theory \cite{rossum}, and the illumination
function on the boundary $g \in H^{1/2}(\partial \Omega)$
satisfies $g \geq 0$ {\it a.e.} on $\partial \Omega$, and
$\dr_\nu$ denotes the normal derivative at $\partial \Omega$.

In diffuse optical tomography, the inverse problem is to
reconstruct the optical absorption distribution $a_*$ from
measurements of the outgoing light intensity on $\partial \Omega$
given by ${\partial_\nu \Phi_{*}} |_{\partial \Omega}$, see
\cite{simon, john}. It is worth mentioning that, in our diffusion
equation model (\ref{Optical eqn}), if $l \neq 0$, then knowing $
\Phi_{*}$ or ${\partial_\nu \Phi_{*}}$ on $\partial \Omega$ is
mathematically the same.

Diffuse optical tomography produces images with poor accuracy and
spatial resolution.  It is known to be ill-posed  due to the fact
that the outgoing light intensities are not very sensitive to
local changes of the optical absorption distribution \cite{simon,
john2,john3,john}.  In
\cite{AmmariBossyLocLaurentetal:preprint2012} we have proposed an
original method for reconstructing the optical absorption
coefficient by using mechanical perturbations of the medium. While
taking optical measurements  the medium is perturbed by a
propagating acoustic wave. Then cross-correlations between the
boundary values of the optical energy density in the medium
changed by the propagation of the acoustic wave and those of the
optical energy density in the unperturbed one are computed.
Finally, under the Born approximation \cite{born}, the use of a
spherical Radon transform inversion yields  a reconstructed image
for $a_*$, which has a resolution of order the width of the wave
front of the acoustic wave propagating in the medium. The Born
approximation linearizes the reconstruction problem. It consists
of assuming that $a_*$ is close to a constant and taking the
background solution of the diffusion equation for constant optical
absorption in place of $\Phi_*$ \cite{john} as the driving optical
energy density at each point in $\Omega$.

The idea of mechanically perturbing the medium has been first
introduced in \cite{AmmariLaurentetal:preprint2012} for
electromagnetic imaging. On the other hand, it is also worth
emphasizing that this approach is different from the imaging by
controlled perturbations \cite{AMMARI-08,
AMMARI-BONNETIER-CAPDEBOSCQ-08,
CAPDEBOSCQ-FEHRENBACH-DEGOURNAY-KAVIAN-09, siap2011, ip12,
GEBAUER-SCHERZER-08, otmar2}, where local changes of the
parameters of the medium are produced by focalizing an ultrasound
beam. Both techniques lead to resolution enhancements. In imaging
by controlled perturbations, the resolution is of order the size
of the focal spot while here it is of the order of the width of
the wave front of the wave propagating in the medium.

This paper aims to generalize the acousto-optic process behind the
Born approximation. We tackle the nonlinear optical reconstruction
problem. We develop a mathematical framework for the
reconstruction problem in the case where the optical absorption
distribution is a perturbation of a piecewise constant function.
We introduce an iterative reconstructing algorithm of
Landweber-type and prove its convergence and stability. For doing
so, we introduce a weak Helmholtz decomposition and interpret in a
weak sense the cross-correlation measurements.

To describe our approach, we employ several notations. Each smooth
component of $a_*$ is called an inclusion. The background of $a_*$
is assumed to be a known positive constant and denoted by $a_0.$
Assume further the knowledge of a lower bound $\underline a$ and
an upper bound $\overline a$ of $a_*$, both of which are positive.
Finally, let $D \Subset \Omega$ be known and such that
\begin{equation} \label{defastar}
    a_* = a_0 \hspace*{.24in} \mbox{in }   \Omega \setminus D.
\end{equation}
We next impose some conditions on the unknown inclusions. Let $k \geq 1$ denote the number of
inclusions and $A_i$ be occupied by the $i$th inclusion. Assume:
\begin{enumerate}
    \item[$I_1$.] for any $i \in \{1, \dots, k\}$, $A_i$ is a smooth subdomain of $\Om,$ $\dr A_i$ is
    connected;
    \item[$I_2$.] for any $j \not = i,$ $\overline{A_i} \cap \overline{A_j} =
    \emptyset$;
    \item[$I_3$.] $\cup_{i = 1}^k \overline A_i \Subset D.$
\end{enumerate}
All of the assumptions above suggest the definition of the class
$(\A) \subset L^\ii(\Om)$, which contains $a_*$.
\begin{de}
    The function $a$ is said to belong to class ($\A$) iff there exist $k \geq 1$,
    $A_1, \cdots, A_k \Subset D$ satisfying $I_1, I_2$ and $I_3$ and $a_1, \cdots,
    a_k \in \C^2(\overline A_i, [\underline a, \overline a])$ such that
    \begin{equation}
    a = \sum_{i = 0}^k a_i\1_{A_i},
    \label{main form}
\end{equation}
where, again, $a_0$ was introduced in (\ref{defastar}), $A_0 =
\Omega \setminus \cup_{i = 1}^k \overline A_i$, and  $\1_{A_i}$
denotes the characteristic function of $A_i$. \label{def 1.1}
%
\end{de}

Our main results in this paper can be summarized  as follows. A
spherical acoustic wave is generated at $y$ outside $\Omega$. Its
propagation inside the medium $\Omega$ changes the optical
absorption distribution. Due to the acoustic wave, any point $x
\in \Omega$ moves to its new position $x+ v_{y, r}^{\eta}(x)$,
where $v_{y, r}^{\eta}$ is defined by (\ref{changes}) with $r$
being the radius of the spherical wave impulsion. By
linearization, the displacement field is approximately $v_{y,
r}^{\eta}$ as the thickness $\eta$ of the acoustic wavefront goes
to zero. Hence, the optical absorption of the medium changed by
the propagation of the acoustic wave is approximately $a_*(x+
v_{y, r}^{\eta})$, up to
 an error of order $\eta$.

Using cross-correlations between the outgoing light intensities in
the medium changed by the propagation of the acoustic wave  and
those of  in the unperturbed one, we get the data  $M_\eta(y, r)$
given by (\ref{data}). In Propositions \ref{approx} and
\ref{approx2}, we show that $M_\eta(y, r)$ converges in the sense
of distributions to $M(y,r)$ as $\eta \rightarrow 0$. We refer to
$M(y,r)$ as the ideal data. Making use of a weak Helmholtz
decomposition, stated in Lemma \ref{lemhelm}, we relate in Theorem
\ref{Theorem Radon psi}  the ideal data to the gradient of
$\Phi_*^2\nabla a_*$. Since $a_*$ is piecewise smooth, $\nabla
a_*$ can be defined only in the  sense of distributions. Technical
arguments and quite delicate estimates are needed in order to
establish the fact that the gradient part of $\Phi_*^2\nabla a_*$
can be obtained from the cross-correlation measurements using the
inverse spherical Radon transform. Based on this, we propose an
optimal control approach for reconstructing the values of $a_*$
inside the inclusions. For doing so, we first detect the support
of $a_*-a_0$ as the support of the gradient part of the data
$\Phi_*^2\nabla a_*$. In fact, Lemma \ref{Lemma Helmholtz
decomposition} shows that the support of the data yields the
support of the inclusions. Their boundaries are detected as the
support of the discontinuities in the data. Proposition
\ref{Theorem 1} provides a Lipschitz stability result for
reconstructing piecewise constant optical absorption. In contrast
with the recent results in \cite{sergio, beretta, beretta2},
Proposition \ref{Theorem 1} uses only one measurement but the
supports of the inclusions are known. Minimizing the discrepancy
functional (\ref{discfunc}) we obtain the background constant
values of the optical absorption inside the inclusions. Next, in
order to recover spatial variations of $a_*$ inside the
inclusions, we minimize the discrepancy between the linear forms
$F[a]$ and $\Delta \psi$ given by (\ref{fa}) and (\ref{defpsif}),
respectively. We prove in Theorem \ref{Lemma F derivative} that
the Fr\'echet derivative of the nonlinear discrepancy functional
is well-defined and establish useful estimates as well. We
introduce an iterative scheme of Landweber-type for minimizing the
discrepancy functional and prove in Theorem \ref{theocv} its
convergence provided that the optical absorption coefficient is in
the set $K$ defined by (\ref{K def}).

\section{Preliminaries}

\subsection{Some basic properties}

We first recall the following results.

\begin{prop}[weak comparison principle \cite{AmmariBossyLocLaurentetal:preprint2012}] \label{prop2.2}
    Let $a \in L^\infty(\Omega)$ be a nonnegative function and assume that $\Phi \in H^1(\Omega)$ satisfies
    \begin{equation}
    \left\{
        \begin{array}{rcll}
            -\Delta \Phi + a \Phi &\geq& 0 &\mbox{in } \Omega,\\
            \nm
            l \partial_{\nu}\Phi + \Phi &\geq& 0 &\mbox{on } \partial \Omega.
        \end{array}
    \right.
     \end{equation}
We have $\Phi \geq 0$ {\it a.e.} in $\Omega$.
\label{weak comparison principle}
\end{prop}

\begin{lemm}[Lemma 4.1 in \cite{AmmariBossyLocLaurentetal:preprint2012}]
    Let $D$ be as in (\ref{defastar}) and assume that $g \in H^{1/2}(\partial \Omega)$ is nonnegative.
  There exist two positive constants
    $\lambda$ and $\Lambda$
    such that for all $a \in (\mathcal{A})$, the solution  $\Phi$ of
\begin{equation}
\left\{\begin{aligned}
    -\La\Phi +a \Phi &=0\ \text{ in } \Om,\\
    \nm
    l \dr_\nu\Phi + \Phi &= g\ \text{ on } \partial\Om,
\end{aligned} \right.
\label{Optical eqn_b}
\end{equation}
satisfies
    \begin{equation}
        \lambda \leq \Phi \leq \Lambda \hspace*{.24in} \mbox{in }
        D.
    \label{2.2}
    \end{equation}
\label{lemma: varphi lambda Lambda}
\end{lemm}

\begin{lemm}[Lemma 4.2 in \cite{AmmariBossyLocLaurentetal:preprint2012}]
Let $T$ be the map that sends $a \in (\mathcal{A})$ into the
unique solution of (\ref{Optical eqn}) with $a$ replacing $a_*$.
Then, $T$ is Fr\'{e}chet differentiable. Its derivative at $a$ is
given by
    \begin{equation}
        DT[a](h) = \varphi,
    \label{2.3}
    \end{equation} for $h \in L^{\infty}(\Omega)$, where $\varphi$ solves
    \begin{equation}
        \left\{
            \begin{array}{rcll}
                -\Delta \varphi + a \varphi &=& -hT[a] &\mbox{\rm in } \Omega,\\
                \nm
                l \partial_{\nu} \varphi + \varphi &=& 0 &\mbox{\rm on } \partial
                \Omega.
            \end{array}
        \right.
    \label{2.4}
    \end{equation} Moreover, $DT[a]$ can be continuously extended to
      $L^2(\Omega)$ by the same formula given in (\ref{2.3}) and (\ref{2.4}) with
\begin{equation}
    \|DT[a]\|_{\mathcal{L}(L^2(\Omega), H^1(\Omega))} \leq
    C\Lambda,
\label{2.5}
\end{equation}
where $\Lambda$ is defined in Lemma \ref{lemma: varphi lambda
Lambda} and $\mathcal{L}(L^2(\Omega), H^1(\Omega))$ is the set of
bounded linear operators from $L^2(\Omega)$ into $H^1(\Omega)$.
\label{Lemma T derivative}
\end{lemm}

The following lemma will be helpful to prove the uniqueness of the
constructed coefficient. We refer to Appendix \ref{appenA} for its
proof.
\begin{lemm} \label{The unique continuation property}

    Let $\Omega'$ be the union of several subdomains of $\Omega$ such that $\Omega \setminus \Omega'$
    is path connected. If $\phi$ is a bounded solution to
    \begin{equation}
        \left\{
            \begin{array}{rcll}
                -\Delta \phi + c \phi &=& 0 &\mbox{in } \Omega \setminus \Omega',\\
                \nm
                l \partial_{\nu}\phi + \phi &=& 0 &\mbox{on } \partial \Omega,
            \end{array}
        \right.
    \label{2.6}
    \end{equation} for some nonnegative constant $c$ and $\partial_\nu \phi \equiv 0$ on $\partial \Omega$,
    then $\phi \equiv 0$
    in $\Omega \setminus \Omega'.$
\label{lemma 2.1}
\end{lemm}

\begin{cor}
    Let $A_0, A_1, \cdots, A_k$ be as in Definition \ref{def 1.1} and let
    $a \in (\A)$ be defined by such sets. Denote by $\varphi_j,$ $j = 1, \cdots, k,$ the solution of
    \[
        \left\{
            \begin{array}{rcll}
                -\Delta \varphi_j + a \varphi_j &=& \1_{A_j}\Phi &\mbox{in } \Omega,\\
                \nm
                l \partial_{\nu} \varphi_j + \varphi_j &=& 0 &\mbox{on } \partial \Omega,
            \end{array}
        \right.
    \] with $\Phi$ being the solution of (\ref{Optical eqn_b}). Then,
    the set $\{\partial_\nu \varphi_j|_{\partial  \Omega}\}$ is linearly independent.
\label{Corollary independent}
\end{cor}
\Proof
   Define
    \[
        \varphi = \sum_{j = 1}^k \alpha_j \varphi_j \hspace*{.24in} \mbox{in }
        \Omega,\] for some $\alpha_1, \cdots, \alpha_k \in \mathbb{R},$
        and assume that $\partial_\nu \varphi = 0$ on $\partial
        \Omega$.  It is
obvious that $\varphi$ is the solution of
    \[
        \left\{
            \begin{array}{rcll}
                -\Delta \varphi + a \varphi &=& \ds \sum_{j = 1}^k \alpha_j\1_{A_j} \Phi &\mbox{in } \Omega,\\
                \nm
               l \partial_{\nu} \varphi + \varphi &=& 0 &\mbox{on } \partial
               \Omega,
            \end{array}
        \right.
    \]  and, hence, satisfies
     (\ref{2.6}) with $c = a_0$ and $\Om' = \cup_{i = 1}^k A_i.$ Thus, by Lemma \ref{lemma 2.1},
     $\varphi \equiv 0$ in $\Omega_0.$ On the other hand, for each $i \in \{1, \cdots, k\},$ $\varphi$ solves
    \[
        \left\{
            \begin{array}{rcll}
                -\Delta \varphi + a_i \varphi &=& \alpha_i\1_{A_i}\Phi &\mbox{in } A_i,\\
                \nm
                \varphi &=& 0 &\mbox{on } \partial A_i.
            \end{array}
        \right.
    \] We can now  apply the strong comparison principle (see, for instance, Lemma 3.1 in \cite{LocSchmitt:die2009}) and
    the Hopf lemma to see that $\partial_{\nu}\phi \not = 0$ on $\partial A_i.$ This contradicts
     to the fact that $\phi \equiv 0$ in $A_0.$ \qed

\subsection{The Helmholtz decomposition in the sense of distributions} \label{Helmholtz decomposition}
    The Helmholtz decomposition plays a crucial role in \cite{AmmariBossyLocLaurentetal:preprint2012}
    when we established a differential coupling system for $a$, where $a$ was supposed to be in $\C^2(\overline \Omega)$.
    Fortunately, when $a$ is no longer smooth but $\Phi^2\nabla a$ belongs to
     $(H^{1}(\Om)^d)^* \subset H^{-1}(\Om)^d$ for all $\Phi \in \C^1(\Omega)$,
      a corresponding Helmholtz decomposition remains true. Note that for all $a \in (\A)$ and $\Phi \in \C^1(\Omega),$
$\Phi^2\nabla a \in (H^{1}(\Om)^d)^*$ in the sense that
\begin{eqnarray}
    \hspace*{-.24in}\langle \Phi^2\nabla a, v \rangle_{(H^{1}(\Om)^d)^*,
    H^{1}(\Om)^d} &=& \langle \Phi^2 \nabla (a - a_0), v \rangle_{(H^{1}(\Om)^d)^*, H^{1}(\Om)^d} \nonumber\\
     &=& -\int_{D}(a - a_0) \Div (\Phi^2 v)\, dx. \label{1.3}
\end{eqnarray}
The domain of the integral above is written as $D$ instead of
$\Omega$ because $a - a_0 = 0$ in $\Om \setminus D$, where $D$ is
introduced in (\ref{defastar}). By the same reason, we do not
require the boundary zero value for the admissible test functions.
The last equation in (\ref{1.3}) suggests that it might be
sufficient to impose $a \in \C^1(\overline A_i)$, instead of
$\C^2(\overline A_i)$, $i = 1, \cdots, k$, as in Definition
\ref{def 1.1}. However, we need the differentiability of $a$ up to
second order in each inclusion for some later regularity and
estimation purposes.

The following result holds.
\begin{lemm} \label{lemhelm}
 For any $U$ in $\hh^d$ there exist  $\psi\in L^2(\Om)$ and $\Psi\in\hh^d$  such that
$$U=\g \psi + \Psi$$
with $\Div \Psi = 0$.
In particular, if $U = \Phi^2\nabla a$ for some $a \in \A$ then $\psi$ is continuous and discontinuous at the
point where $a$ is, respectively. 
\label{Lemma Helmholtz decomposition}
\end{lemm}

\Proof Letting $U = (U_1, \cdots, U_d) \in H^{-1}(\Omega)^d,$ we
denote by $u = (u_1, \cdots, u_d)$ the solution of
\begin{equation}
    \left\{
        \begin{array}{rcll}
            -\Delta u &=& U &\mbox{in } \Omega,\\
            u &=& 0 &\mbox{on } \partial
\Omega.
        \end{array}
    \right.
\label{Riezt}
\end{equation}
The vector $u \in H^1_0(\Omega)^d$ is actually the Riesz
representation of $U$ in $H^1_0(\Omega)^d.$ Applying the classical
Helmholtz decomposition for $u$ (see, for instance, \cite{galdi}),
we can find $f \in H^1(\Omega)$ and $G \in H(\mbox{curl},
\Omega):=\{ w \in L^2(\Omega)^d: \, \nabla \times w \in
L^2(\Omega)^d\}$ such that
\begin{equation}
    u = \nabla f + \nabla \times G.
\label{1}
\end{equation}
Here, $\Div G = 0$ inside $\Omega$ and \begin{equation}
\label{gtimesnu} G\times \nu =0 \quad \mbox{on }
\partial \Omega.\end{equation} Moreover, $f$ is a solution of
\begin{equation}
\left\{
\begin{array}{rcll}
    \Delta f &=& \Div u &\mbox{in } \Omega,\\
    \nm
        \partial_\nu f &=& 0 &\mbox{on } \partial \Omega.
\end{array}
    \right.
\label{2.10}
\end{equation}
Since $u$ belongs to $H^1(\Omega)^d$, $\Div u \in L^2(\Omega)$. By
standard regularity results, we see that $f \in H^2(\Omega).$


In view of (\ref{Riezt}), taking the Laplacian of (\ref{1}) yields
\[
    U = \nabla \psi + \Psi,
\] in the sense of distributions, where $\psi = \Delta f \in L^2(\Omega)$ and $\Psi$ is divergence free.


We next prove the second statement of the lemma in which $U =
\Phi^2\nabla a$ for some $a \in (\mathcal{A})$. The main tools we
use here are the $H^2$- and $\C^1$-regularity results. Fix $j \in
\{1, \cdots, d\}$ and $i \in \{0, \cdots, k\}$. Denote by $u_j$
the $j$th component of the vector $u$, defined in (\ref{Riezt}).
Since $u_j \in H^1_0(\Omega)$, it belongs to $H^1(A_i)$. The
function $u_j$ solves
\begin{equation}
    -\Delta u_j = \Phi^2\partial_{x_j}a,
\label{2.11}
\end{equation}
in $A_i$. Applying Theorem 8.8 in \cite{GilbargTrudinger:1977}, we
see that $u_j$ is in $H^2(A_i')$ for all $A_i' \Subset A_i$.
Hence, differentiating (\ref{2.11}) gives
\[
-\Delta \partial_{x_l} u_j = \partial_{x_l}(\Phi^2\partial_{x_j}a)
\] in $A_i'$ for all $l = 1, \cdots, d.$ Since $\partial_{x_l} u_j \in H^1(A_i')$
and $\partial_{x_l}(\Phi^2\partial_{x_j}a) \in L^2(A_i')$, we can
apply the $\C^1$-regularity result in
\cite{LadyzhenskayaUraltseva:sv1985} to see that $\partial_{x_l}
u_j$ is in $\C^1(A_i'')$ for all $A_i'' \Subset A_i'.$ This
implies $u_j \in \C^2(A_i)$. Considering the differential equation
in (\ref{2.10}) in each inclusion and following the same
regularity process, we see that $f \in \C^2(A_i)$. Hence $\psi =
\Delta f$ is continuous in $A_i,$ which is also the set of
continuous points of $a$. On the other hand, since $U =
\Phi^2\nabla a$ involves Dirac distributions supported in $\cup_i
\partial A_i$, $\nabla \cdot u$ is not continuous across $\cup_i
\partial A_i$, so are $f$ and $\psi = \Delta f$. \qed


\section{The set of data}
In this section, we describe the set of data obtained by the
acousto-optic process introduced in
\cite{AmmariLaurentetal:preprint2012}. The basic idea in order to
achieve a resolution enhancement in imaging the optical absorption
distribution is a s follows. We generate a spherical acoustic wave
inside the medium. The propagation of the acoustic wave changes
 the absorption parameter of the medium. During the propagation of
 the wave we measure the light intensity on  $\partial \Omega$. The
 aim is now to reconstruct the optical absorption
 coefficient from such set of measurements.

Let $a \in (\mathcal{A})$ represent the true coefficient $a_*$.
Let $S^{d - 1}$ be the unit sphere in $\R^d$. Let $\mu
>0$ and let $S_{\mu} = \mu S^{d - 1}$, the sphere of radius $\mu$
and center $0$,  be such that $\Omega$ stays inside $S_{\mu}.$ We
perturb the optical domain $\Omega$ by spherical acoustic waves
generated at point sources $y \in S_{\mu}$. Let $r \in [r_0, R]$
be the radius of the spherical wave impulsion, where $r_0$ and $R$
are the minimum and maximum radii so that the spherical waves
generated at point sources on $S_{\mu}$ can intersect $\Om$. Let
$\eta \ll 1$ be the acoustic impulsion typical length representing
the thickness of the wavefront. Let the position function $P$ be
defined by
\[
    P : x \mapsto x + v_{y, r}^{\eta} (x), \quad x \in
    \Omega,
\] where
\begin{equation} \label{changes}
    v_{y, r}^{\eta} (x)=
    \eta\frac{r_0}{r}w\left(\frac{r-|x-y|}{\eta}\right)\frac{x-y}{|x-y|},
\end{equation} and $w$ is a smooth function supported on $[-1,1]$ with $\n w
\ii=1$. Here, $\n \, \ii$ denotes $\n \, {L^\infty(]-1,1[)}$.

In \cite{AmmariLaurentetal:preprint2012}, we have shown that the
displacement function  at the point $x$ caused by the short
diverging spherical acoustic wave generated at $y$ is given by
\begin{equation} \label{changes2}
     u_{y,r}^\eta(x)=  P^{-1}(x) - x, \quad x \in \Omega.
\end{equation}

 Let $C$ be the cylinder $S_{\mu} \times [r_0, R]$. For each $(y, r) \in C$, $a_{u_{y,r}^\eta}(x)$
 denotes $a(x+u_{y,r}^\eta(x))$ and $\Phi_{u_{y,r}^\eta}$ is the optical energy density in the
 displaced medium, which satisfies
\begin{equation}
\left\{\begin{aligned}
-\La\Phi_{u_{y,r}^\eta}+a_{u_{y,r}^\eta}\Phi_{u_{y,r}^\eta} &=0\ \text{ in } \Om,\\
l \dr_\nu\Phi_{u_{y,r}^\eta}+ \Phi_{u_{y,r}^\eta} &= g\ \text{ on
}
\partial\Om.
\end{aligned} \right.
\end{equation}
Physically, the outgoing light intensities $\partial_\nu \Phi
|_{\partial \Om}$ and $\partial_\nu \Phi_{u_{y,r}^\eta}|_{\partial
\Om}$ are measured. We are thus able to assume the knowledge of
the cross-correlation measurements:
\begin{equation} \label{cross}
    \frac{1}{\eta^2}\int_{\partial \Omega}g (\partial_\nu \Phi - \partial_\nu \Phi_{u_{y,r}^\eta})d\sigma, \hspace*{.24in}
    y \in S_{\mu}, r > 0.
\end{equation}
Integration by parts shows that the quantity above is equal to
\begin{equation} \label{data}
    M_\eta(y, r) = \frac{1}{\eta^2}\int_\Om (a_{u_{y, r}^\eta} - a)\Phi\Phi_{u_{y, r}^\eta}dx,
\end{equation}
which is considered as our set of data. Here, the coefficient
${1}/{\eta^2}$ is put in front of the integral because both $\Phi$
and $\Phi_{u_{y, r}^\eta}$ are bounded (Lemma \ref{lemma: varphi
lambda Lambda}) and
\begin{equation}
    \|a_{u_{y, r}^\eta} - a\|_{L^1(\Om)} = O(\eta^2) \hspace*{.24in } \mbox{as } \eta \rightarrow
    0^+,
\label{3.5}
\end{equation}
provided that the following technical condition, named as ($\mathcal{H}$), is imposed: there exists $\delta > 0$
such that for all $x \in \partial A_i \cap \Sigma_\eta(y, r)$, either
\begin{itemize}
    \item[$H_1$:] the angle formed by the ray $x - y$ and the normal outward vector of $A_i$ at $x$ is
    greater than $\delta$; or,
    \item[$H_2$:] the curvature of $\partial A_i$ is different to that of the circle
    or sphere $\{z \in \R^d: |z - y| = |x - y|\}$ at $x$ if the angle above is smaller than $\delta$. \
\end{itemize} Here,
\[
    \Sigma_\eta(y, r) = \{z \in \R^d: r - \eta < |z - y| < r + \eta\}.
\]
In fact, this condition guarantees that
\begin{equation}
    |A_i \La P^{-1}(A_i)| + |A_i \La P(A_i)|\leq O(\eta^2).
\label{2.4'}
\end{equation}
Denote
\begin{equation}
    V_\e(S)=\{x\in\R^d,\ \exists y \in S,\ |x-y|<\e\},
    \label{2.4444}
\end{equation} for any smooth surface $S$ of $\R^d$, and $\e > 0.$
Since $S$ is smooth, the volume of $V_\e(S)$ is given by
$$V_\e(S)= 2\s(S)\e+ O(\e^2).$$
Fix $(y, \eta) \in C$ and write
\begin{eqnarray}
\n{a_{u_{y, r}^\eta} - a}{L^1(\Om)} &=& \sum_{i=1}^n\int_{A_i\cup P^{-1}(A_i)}|a_{u_{y, r}^\eta} - a|dx \nonumber \\
&=& \sum_{i=1}^n \int_{A_i\cap P^{-1}(A_i)}|a_{u_{y, r}^\eta} - a|dx + \int_{A_i \La P^{-1}(A_i)}|a_{u_{y, r}^\eta}
- a|dx \label{3.6}.
\end{eqnarray}
As $u_{y, r}^{\eta}$ is supported on $\Sigma_\eta(y, r)$ and $\n{u_{y, r}^{\eta}}{\ii} = \eta$,
\begin{eqnarray*}
\int_{A_i\cap P^{-1}(A_i)}|a_{u_{y, r}^{\eta}}-a|dx &=& \int_{\Sigma_\eta\cap A_i\cap P^{-1}(A_i)}|a_{u_{y, r}^{\eta}} - a|dx
\\&\leq& \eta\n {\g a_i}{L^\ii(A_i)}|\Sigma_\eta|\\
&\leq& \n {\g a_i}{L^\ii(A_i)}\s(S(0,R))\eta^2,
\end{eqnarray*}
where $\s(S(0,R))$ is the surface measure of the sphere of center
$O$ and radius $R$.  The second integral in (\ref{3.6}) is bounded
by $O(\eta^2)$ because of (\ref{2.4'}) and the boundedness of $a$.

\section{The behavior of $M_{\eta}$ as $\eta$ approaches $0^+$ and the ideal measurements}

Consider the open cylinder $C:=S_\mu\times (0, R)$ with its
classical product topology.

 The construction of
\begin{eqnarray*}
    M_{\eta} : C &\rightarrow& \R \\
                 (y, r) &\mapsto& \frac{1}{\eta^2}\int_{\Omega} (a_{u^{\eta}_{y, r}} - a)\Phi \Phi_{u^{\eta}_{y, r}} dx,
\end{eqnarray*}
has been described in this previous section. The knowledge of this
function is obtained from those of $g,$ $\partial_\nu \Phi$ and
$\partial_\nu \Phi_{u^{\eta}_{y, r}}$ on $\partial \Omega$. In
this  section, we study the limit of $M_{\eta}$ as $\eta
\rightarrow 0^+$. This, together with a weak version of Helmholtz
decomposition and the spherical Radon transform, will help us to
detect all inclusions.

\begin{lemm}
For any $\eta>0$, $M_\eta$ is a continuous map on $C$.
\label{Lemma 3.1}
\end{lemm}
\Proof It is sufficient to consider only the case $r > r_0$
because $M_{\eta}(y, r) = 0$ for all $r \leq r_0$ and $y \in
S_{\mu}.$ Fix $(y, r) \in S_\mu \times (r_0, R)$ and let $\{(y_n,
r_n)\}_{n \geq 1} \subset S_\mu \times (r_0, R)$ converge to $(y,
r).$ Noting that $a_{u_{y, r}}^{\eta}$ is continuous except on the
zero measured set
\[
    \{x + u_{y, r}^\eta(x): x \in \cup_{i = 1}^n \partial A_i\},
\] we have
\[
    a(x + u_{y_n, r_n}^\eta(x)) \rightarrow a(x + u^\eta_{y, r}(x))
\] {\it a.e.} in $\Omega.$ On the other hand, since $a$ is uniformly
bounded, so is \[|a(x + u^\eta_{y_n, r_n}(x)) - a(x + u^\eta_{y, r}(x))|^2.\] It follows by
the Lebesgue dominated convergence theorem that
\[
    a_{u^{\eta}_{y_n, r_n}} \rightarrow a_{u^\eta_{y, r}} \hspace*{.24in} \mbox{in } L^2(\Omega)
\]  as $n \rightarrow \ii.$
 This implies
\[\Phi_{u^{\eta}_{y_n, r_n}} \rightarrow \Phi_{u^{\eta}_{y, r}}\]
in both $H^1(\Omega)$ and $L^4(\Omega).$ Note that the $L^4$
convergence above is valid because $d$ is either $2$ or $3$. A
direct calculation yields
\begin{eqnarray*}
    &&|\eta^2(M_{\eta}(y_n, r_n) - M_{\eta}(y, r))|\\
    && \hspace*{.24in} = \left|\int_{\Omega} [(a_{u^{\eta}_{y_n, r_n}} - a)\Phi\Phi_{u^{\eta}_{y_n, r_n}} - (a_{u^{\eta}_{y, r}} - a)\Phi\Phi_{u^{\eta}_{y, r}}]dx\right|\\
    && \hspace*{.24in} \leq \int_{\Omega} |a_{u^{\eta}_{y_n, r_n}} - a|\Phi|\Phi_{u^{\eta}_{y_n, r_n}} - \Phi_{u^{\eta}_{y, r}}|dx + \int_{\Omega} |a_{u^{\eta}_{y_n, r_n}} - a_{u^{\eta}_{y, r}}|\Phi\Phi_{u^{\eta}_{y, r}}dx,\\
    && \hspace*{.24in} \leq 2\overline a \|\Phi\|_{L^4(\Omega)}\|\Phi_{u^{\eta}_{y_n, r_n}} - \Phi_{u^{\eta}_{y, r}}\|_{L^4(\Omega)} \\
    &&\hspace*{2in}+ \|a_{u^{\eta}_{y_n, r_n}} -
     a_{u^{\eta}_{y, r}}\|_{L^2(\Omega)}\|\Phi\|_{L^4(\Omega)}\|\|\Phi_{u^{\eta}_{y, r}}\|_{L^4(\Omega)}\|.
\end{eqnarray*}
The lemma follows. \qed

Lemma \ref{Lemma 3.1} guarantees that $M_\eta$ is measurable. In
the case that $a$ is smooth, which has been studied in
\cite{AmmariLaurentetal:preprint2012,AmmariBossyLocLaurentetal:preprint2012},
$M_\eta(y,r)\approx \int_\Om\g a \cdot u_{y,r}^\eta\Phi^2$ when
$\eta$ is small. However, when $a$ is piecewise smooth, we need to
establish a similar approximation in the weak sense. The following
proposition holds. We refer to Appendix \ref{appenB} for its
proof.

\begin{prop} \label{approx}
Let $C=S_\mu\times (0, R)$. For any $0 < \eta \ll 1$, define the
continuous function
\begin{equation}
 \ti M_\eta(y,r)=\frac 1 {\eta^2}\int_\Om (a-a_0)\Div(\Phi^2v_{y,r}^\eta)dx, \hspace*{.24in} (y, r) \in
 C,
 \label{tilde M}
\end{equation}
where  $\Phi$ is the solution of (\ref{Optical eqn}) with $a$
replacing $a_*$. Assume ($\mathcal{H}$) holds and, consequently,
(\ref{2.4'}) is valid. Then there exists $c>0,$ independent of
$(y, r)$, such that
\begin{equation}
\left|M_\eta(y,r) - \ti M_\eta(y,r)\right|\leq c\eta,\ \ \forall(y,r)\in C.
\end{equation}
\end{prop}

It follows from Proposition \ref{approx} that for each $(y, r) \in
C,$
\begin{equation}
    \lim_{\eta \rightarrow 0^+} M_{\eta}(y, r) = \lim_{\eta \rightarrow 0^+}
    \ti M_\eta(y,r)  := M_{y, r}.
\label{3.4}
\end{equation}
We cannot expect that $M$ is a smooth function on $C$ because
${u^\eta_{y, r}}/{\eta^2}$, and hence ${v^\eta_{y, r}}/{\eta^2}$,
converges to a distribution supported on the circle (or sphere)
$S(y, r) = \{z : |z - y| = r\}$. The limit in (\ref{3.4}) is
understood as follows.

Let
\[
    G(C) = \left\{f \in L^2(C): \partial_r f \in L^2(C)\right\},
\] be a Hilbert space, endowed with the norm
\[
    \|\cdot\|_{G(C)} = \|\cdot\|_{L^2(C)} + \|\partial_r \cdot\|_{L^2(C)}.
\] Let $\gamma$ be the (continuous) trace operator from $C$ to $S_\mu \times \{0, R\}$ and denote
\[
    G_0(C) = \gamma^{-1}(0) = \{f \in G(C): \gamma(f) = 0\}, \hspace*{.24in} G^{-1}(C) = G_0(C)^*.
\] We have the following relations
\[
    H^1_0(C) \subset G_0(C) \subset L^2(C), \hspace*{.24in} L^2(C) \subset G^{-1}(C) \subset H^{-1}(C).
\]

Let $\n \, 1$ denote $\n \, {L^1(]-1,1[)}$. The following is the
main result of this section. It is a direct consequence of
Proposition \ref{approx}.
\begin{prop} \label{approx2}
The function ${M_\eta}$ converges to the ideal measurements $M$ in
$G^{-1}(C)$ as $\eta \rightarrow 0^+$ with
\begin{eqnarray}
    &&\hspace*{-.24in}\langle M,\ph \rangle_{(\mathcal{C}_0^\infty(C))^*, \mathcal{C}_0^\infty(C)} \nonumber\\
    &&\hspace*{.12in} = - r_0\n w 1\int_{S_\mu}\int_0^R\int_{S^{d-1} \cap \Omega_{y, r}}a(y+r\xi)
    \frac{\dr}{\dr r}\left(r^{d-2}\Phi^2(y+r\xi)\ph(y,r)\right)d\xi drdy , \label{4.5}
\end{eqnarray}
where
\[
    \Omega_{y, r} = \left\{\frac{x - y}{r}: x \in \Omega \right\}.
\]
\end{prop}
\Proof For any $\ph$ in $G_0(C)$, we have
\begin{equation}\begin{aligned}\nonumber
\left<\ti M_\eta,\ph\right> &=-\int_{y\in S_\mu}\int_{r=0}^R\int_\Om(a-a_0)(x)\Div_x\left(\Phi^2(x)\frac{v_{y,r}^\eta(x)}{\eta^2}\ph(y,r)\right)dxdrdy\\
&=-\int_{y\in S_\mu}\int_\Om(a-a_0)(x)\Div_x\left(\Phi^2(x)\int_{r=0}^R\frac{v_{y,r}^\eta(x)}{\eta^2}\ph(y,r)dr\right)dxdy.\\
\end{aligned}\end{equation}
Then by the change of variables $x=y+\rho\xi$ we get
$$\frac{v_{y,r}^\eta(x)}{\eta^2}=\frac{r_0}{r\eta}w(\frac{\rho-r}{\eta})\xi.$$
Hence we can write
\begin{equation}\begin{aligned}\nonumber
&\left<\ti M_\eta,\ph\right> = -\\
&\int_{S_\mu}\int_{S^{d-1}}\int_{\rho=0}^R(a-a_0)(y+\rho\xi)\frac{\dr}{\dr \rho}\left(\rho^{d-1}\Phi^2(y+\rho\xi)\int_{r=0}^R\frac{r_0}{r\eta}w\left(\frac{\rho-r}{\eta}\right)\ph(y,r)dr\right)d\rho d\xi dy\\
\end{aligned}\end{equation}
Since $$\frac{1}{\eta}w\left(\frac{\rho-r}{\eta}\right)\
\overset{\eta\to 0}\fg\ \n w 1\del_\rho,$$ we deduce that
$$\int_{r=0}^R\frac{r_0}{r\eta}w\left(\frac{\rho-r}{\eta}\right)\ph(y,r)dr\
\overset{\eta\to 0}\fg\ \frac{\n w 1 r_0} \rho\ph(y,\rho),$$ and
then
\begin{equation}\begin{aligned}\nonumber
&\left<\ti M_\eta,\ph\right> \ \overset{\eta\to 0}\fg\ \\
&-\n w 1 r_0\int_{S_\mu}\int_{\rho=0}^R\int_{S^{d-1}}
a(y+\rho\xi)\frac{\dr}{\dr
\rho}\left(\rho^{d-2}\Phi^2(y+\rho\xi)  \ph(y,\rho)\right)d\rho d\xi dy,\\
\end{aligned}\end{equation}
as desired. \qed

\section{Detecting the inclusions} \label{sec: Detecting the inclusions}
%
Using the fact that $\Phi^2\g a \in (H^1(\Om)^d)^* \subset \hh^d$,
we can employ Lemma \ref{Lemma Helmholtz decomposition} to write
that
\begin{equation}
\Phi^2\g a=\g \psi + \Psi,  \label{5.1}
\end{equation}
where $\Psi$ is a divergence free field and $\psi \in L^2(\Om)$.
Since that both $\Phi^2\g a$ and $\g \psi$ are in
$(H^1(\Om)^d)^*,$ so is $\Psi$. Moreover, it follows from the
usual integration by parts formula and the boundary condition
(\ref{gtimesnu}) that
\begin{equation}
    \langle \Psi, \nabla v \rangle = 0, \hspace*{.24in} \forall\; v \in \C^{\ii}(\overline \Omega).
    \label{5.2}
\end{equation}

For a distribution $f \in (\mathcal{C}_0^\infty(C))^*$, we define
its spherical Radon transform $\Ra[f]$ in the sense of
distributions by
$$ \langle \Ra[f],\ph \rangle_{(\mathcal{C}_0^\infty(C))^*,
\mathcal{C}_0^\infty(C)}  = \langle f, \Ra^*[\ph]
\rangle_{(\mathcal{C}_0^\infty(C))^*, \mathcal{C}_0^\infty(C)},$$
where $$\Ra^*[\ph](x) = \int_C  \ph(y, |x-y|) dy \quad \mbox{for }
\ph \in \mathcal{C}_0^\infty(C).$$ We have the following result.
\begin{theo}
The spherical Radon transform $\Ra[\psi]$ of $\psi$ satisfies the
equation
\begin{equation} \label{formulainvert}
M= r_0\n w 1 r^{d-2} \frac{\dr \Ra[\psi]} {\dr r}
\end{equation}
in the sense of distributions. \label{Theorem Radon psi}
\end{theo}

\Proof Let $\ph\in \C^{\ii}_0(C)$, and for a fixed $y\in S_\mu$ we
define
    $$F_y(x) = \ph(y,|x-y|)\frac{x-y}{|x-y|^2} \hspace*{.24in} x \in \Omega.$$
For any $y \in S_\mu$, the vector $F_y$ is in $H^1(\Om)^d$ because
$|x - y| \geq r_0$. Equation (\ref{5.1}) yields
$$\langle \Phi^2\g a,F_y \rangle = \langle \g\psi,F_y \rangle + \langle \Psi,F_y\rangle$$
and since $F_y$ is the gradient of the function given by
$$x\longmapsto\int_0^ {|x-y|}\frac{\ph(y,\rho)}{\rho} d\rho,$$ it
follows from (\ref{5.2}) that $\langle \Psi ,F_y \rangle = 0$.
Here,  $\langle \, ,\, \rangle$ denotes the duality pair between
$H^1(\Om)^d$ and $(H^1(\Om)^d)^*$. Therefore,
\begin{equation}
    \langle \Phi^2\g a,F_y \rangle = \langle \g\psi,F_y \rangle .
\label{5.4}
\end{equation}
A simple calculation shows
\begin{eqnarray*}
    \langle \Phi^2\g a,F_y \rangle &=& \int_\Om a\Div(\Phi^2F_y)dx\\
                                    &=& \int_0^R\int_{S^{d-1} \cap \Omega_{y,r}}\left[a\Div(\Phi^2F_y)\right](y+r\xi)r^{d-1}d\xi
                                    dr,
\end{eqnarray*}
and hence,
\begin{equation} \label{th45}
\langle \Phi^2\g a,F_y \rangle = \int_0^R\int_{S^{d-1}
 \cap \Omega_{y,r}}a(y+r\xi)\frac{\dr}{\dr
r}\left[\Phi^2(y+r\xi)\ph(y,r)r^{d-2}\right]d\xi dr.
\end{equation}
Combining  (\ref{4.5}), (\ref{5.4}), and  (\ref{th45}) implies
\begin{eqnarray*}
    \langle M, \varphi \rangle
    &=& \n{w}{L^1(\Omega)}\int_{S_{\mu}}  \langle \g\psi,F_y \rangle dy\\
    &=& -\n{w}{L^1(\Omega)}\int_{S_{\mu}}  \int_\Om \psi \Div(F_y)dx dy \\
    &=& -\n{w}{L^1(\Omega)}\int_{S_{\mu}}  \int_0^R\int_{S^{d-1}  \cap \Omega_{y,r}}\left[\psi\Div(F_y)\right](y+r\xi)r^{d-1}d\xi dr dy \\
    &=& -\n{w}{L^1(\Omega)}\int_{S_{\mu}} \int_0^R\int_{S^{d-1}  \cap \Omega_{y,r}}\psi(y+r\xi)\ \frac{\dr}{\dr r}\left[\ph(y,r)r^{d-2}\right]d\xi dr dy\\
    &=& - r_0\n w {L^1(\Omega)} \int_{S_\mu}\int_0^R \Ra[\psi](y,r)\frac{\dr}{\dr r}\left[\ph(y,r)r^{d-2}\right] dr dy\\
    &=& r_0\n w {L^1(\Omega)} \langle r^{d-2} \frac{\dr \Ra[\psi]} {\dr
    r},\ph\rangle,
\end{eqnarray*}
and the proof is complete.\qed

\begin{rem}
    Theorem \ref{Theorem Radon psi} provides the knowledge of the derivative of the spherical Radon transform
    of $\psi$ (see Appendix
    \ref{A primitive construction in the sense of distributions} for
    the reconstruction of $\Ra[\psi]$ from its derivative).
Note that the function $\psi$ itself can be reconstructed in a
stable way from $\Ra[\psi]$ using an inversion (filtered)
retroprojection formula for the spherical Radon transform.
    From this,
    all inclusions are detected by the second statement in Lemma \ref{Lemma Helmholtz
    decomposition},
    noticing that $\partial A_i$ is the set of discontinuous points of $\psi.$
\end{rem}

\section{A reconstruction algorithm of the true coefficient}

With all inclusions $A_1, A_2, \cdots, A_k$ in hand, we are able
to find an initial guess for $a_*$ using the unique continuation
property (Lemma \ref{The unique continuation property}) and then
employ a Landweber type iteration to reconstruct $a_*$. As an
initial guess, we reconstruct constant values inside each
inclusion by minimizing the discrepancy between computed and
measured boundary data. We prove a Lipschitz stability result for
the reconstruction of the optical absorption coefficient in the
class of piecewise constant distributions provided that the
support of the inclusions is known.

\subsection{The data of boundary measurements and an initial guess}
\label{sec:simple}

Define
\[
    \mathcal{S} = \left\{\sum_{i = 0}^k \alpha_i \1_{A_i}: \alpha_0 = a_0 \mbox{ and } \alpha_1, \cdots, \alpha_k \in [\underline a, \overline a]\right\}.
\]

Let $a_1$ and $a_2$ be in $\mathcal S$. Their difference can be written as
\[
    a_2 - a_1 = \sum_{i = 1}^k h_i\1_{A_i},
\] for some $h = (h_1, \cdots, h_k) \in B = [\underline a - \overline a, \overline a - \underline a]^k$.
Note that $B$ can be considered as a closed ball of $\mathbb{R}^k$
with respect to the $\infty-$norm of $\mathbb{R}^k$ given by
\[
    |h| = \max\{|h_1|, \cdots, |h_k|\}.
\]  The compactness of $B$ plays an important role in our
analysis. Suppose that $l\neq 0$.  Denote by $\Phi_1$ and $\Phi_2$
the optical energy density functions that correspond to $a_1$ and
$a_2$. The function $\phi = \Phi_1 - \Phi_2$ solves
\begin{equation}
    \left\{
        \begin{array}{rcll}
            -\Delta \phi + a_1\phi &=& \ds \sum_{i = 1}^k h_i\1_{A_i}\Phi_2 &\mbox{in } \Omega,\\
            \nm
            l \partial_{\nu} \phi + \phi &=& 0 &\mbox{on } \partial \Omega.
        \end{array}
    \right.
\label{4.1}
\end{equation}
Using $\phi$ as the test function in the variational form of (\ref{4.1}), we see that
\[
    \int_{\Omega} (|\nabla \phi|^2 + \underline a \phi^2)dx +
    {l} \int_{\partial \Omega} (\partial_\nu \phi)^2 \, d\sigma \leq |h|\lambda
    \int_{\Omega}|\phi|dx,
\]
where $\lambda$ is defined in Lemma \ref{lemma: varphi lambda
Lambda}. This implies
\[
    \|\partial_\nu \phi\|_{L^2(\partial \Omega)} \leq C |h|
\]
and, therefore, the continuity of the map $h \mapsto
\partial_\nu \phi|_{\partial \Omega}$. Since the map $h \in
\partial_{\mathbb{R}^k}B \mapsto \|\partial_\nu\phi\|_{L^2(\partial \Omega)}$
is continuous and nonzero (due to Corollary \ref{Corollary
independent}), we can employ the compactness of
$\partial_{\mathbb{R}^k}B$ in $\mathbb{R}^k$ to see that
\[
    c(a_1) = \min_{h \in \partial_{\mathbb{R}^k}B }  \|\partial_\nu \phi\|_{L^2(\partial \Omega)} > 0.
\]
Identifying $\mathcal{S}$ with a compact subset of $\mathbb{R}^k$, we can conclude that
\[
    c = \inf_{a_1 \in \mathcal S} c(a_1) > 0.
\]
Properly scaling the inequality
\[
    \|\partial_\nu \phi\|_{L^2(\partial \Omega)} \geq c
\] for all $h \in \partial_{\mathbb{R}^k}B$, we arrive at the
following Lipschitz stability result using only one measurement.
Note here that the support of the inclusions is known and only the
value of the optical absorption coefficient inside each inclusion
is to determine.
\begin{prop}
    There exists $c >0$ such that for all $a_1, a_2 \in \mathcal{Q},$
        \begin{equation}
            \|\partial_\nu \Phi_1 - \partial_\nu \Phi_2\|_{L^2(\partial \Omega)} \geq c\|a_1 -
            a_2\|_{L^{\infty}(\Omega)},
        \label{4.2}
        \end{equation} where $\Phi_1$ and $\Phi_2$ are the solutions of (\ref{Optical eqn}) with $a_*$
        replaced by $a_1$ and $a_2$, respectively.
\label{Theorem 1}
\end{prop}

\begin{rem}
    Inequality (\ref{4.2}) guarantees the uniqueness
    of the reconstruction for $a_* \in \mathcal S$ if $\partial_\nu \Phi_*|_{\partial  \Omega}$
    is considered as the data given. It, moreover, implies the stability in
    the sense that small noise does not cause large error.
\end{rem}

Proposition \ref{Theorem 1} suggests us to minimize the quadratic
misfit functional:
\begin{equation} \label{discfunc}
    J(a) = \frac{1}{2} \|\partial_\nu \Phi - \partial_\nu \Phi_*\|^2_{L^2(\partial \Omega)},
\end{equation} where $a$ varies in $\mathcal S$ and $\Phi_*$ is the true
optical energy density. This is possible since $\mathcal S$ is
identical with a compact subset of $\mathbb{R}^k$. By (\ref{4.2}),
the function $a_I = \mbox{argmin }J$ is close to $a_*$ provided
that $a_*$ is a perturbation of a constant on each inclusion
$A_i$. Therefore, $a_I$ can be considered as the background
constant optical absorption distribution in the inclusions. For
simplicity, we propose the following exhaustion method: for each
fine partition $P$ of the interval $[\underline a, \overline a]$,
try all values of $\alpha_i$ such that $\alpha_i$ equals each
element of $P$, and finally choose the $k-$tuple $(\alpha_1,
\cdots, \alpha_k)$ that gives the smallest $\|\partial_\nu \Phi -
\partial_\nu \Phi_*\|_{L^{2}(\partial \Omega)}.$

\subsection{Internal data map and its differentiability}

 Define the set that $a_* = (a_1^*, \cdots, a_k^*)$, identifying with true optical absorption coefficient
 $a_*$  of the form (\ref{main form}), belongs to
\begin{equation}
    K := \{a \in \prod_{j = 1}^k W_0^{1, 4}(A_j): \underline a \leq  a_i \leq \overline
     a \mbox{ and } \|\nabla a_i\|_{L^4(A_j)} \leq \theta, i = 1, \cdots,
     k\},
\label{K def}
\end{equation}
where $\theta$ will be determined later in (\ref{Theta condition}). It is obvious that $K$ is closed
and convex in $H$ where $H = \prod_{j = 1}^k H^1_0(A_j)$ is a Hilbert space with the usual inner product
\[
    \langle u, v \rangle_H = \sum_{i = 1}^k \int_{A_j} \nabla u_j \cdot \nabla v_j dx
\] for all $u = (u_1, \cdots, u_k)$ and $v = (v_1, \cdots, v_k)$ in $H.$

Now, let the map $F: K \rightarrow H^*$ be defined as follows. For
all $(a_1, \cdots, a_k) \in K,$ let
\begin{equation}
    a = \sum_{i = 0}^ka_i \1_{A_i},
\label{q in terms of a}
\end{equation}
and
\begin{equation} \label{fa}
    F[a](v) = \sum_{j = 1}^k\int_{A_j}T[a]^2\nabla a_j \cdot \nabla v \hspace*{.24in} \mbox{for all } v \in
    H,
\end{equation}
where $T[a]$ was defined in Lemma \ref{Lemma T derivative}. We
call $F$ the \textit{internal data} map.

\begin{theo}
    The map $F$ is Fr\'echet differentiable in $K$ and
\begin{equation}
    DF[a](h, v) = \sum_{i = 1}^k \int_{A_i} (2T[a]DT[a](h) \nabla a_i + T[a]^2\nabla h_i) \nabla v_j dx
\label{F derivative}
\end{equation}
for all $a = (a_1, \cdots, a_k) \in K,$ $h = (h_1, \cdots, h_k) \in \prod_{j = 1}^k W^{1, 4}_0(A_j)
\cap L^{\infty}(A_j)$ and $v = (v_1, \cdots, v_k) \in H.$
Assume further
\begin{equation}
    0 < \theta < \frac{C_{\Om'}\lambda^2}{\Lambda^2},
\label{Theta condition}
\end{equation} where $\Om' = \cup_{j = 1}^k A_j$ and $C_{\Om'}$ is the norm of the
embedding map of $H^1(\Om')$ into $L^4(\Om'),$ multiplied with the
constant in (\ref{2.5}). Then, $DF[a]$ is well-defined on $H$ and
there exists a positive constant $C$ such that for all $h \in H$,
\begin{equation}
    \|DF[a](h)\|_{H^*} \geq C\|h\|_{H}.
\label{5.4.2}
\end{equation}
\label{Lemma F derivative}
\end{theo}
\begin{rem}
    The term $DT[a](h)$ in (\ref{F derivative}) is understood as $DT[a]$ acting on the
     function that is equal to $0$ in $A_0$ and to $h_j$ in $A_j$, $j = 1, \cdots, k$.
\end{rem}

\noindent {\it Proof of Theorem \ref{Lemma F derivative}}. The
Fr\'echet differentiability of $F$ and the expression (\ref{F
derivative}) of $DF$ can be deduced from Lemma \ref{Lemma T
derivative} and the standard rules in differentiation. We only
prove (\ref{5.4.2}). In fact, for all $h \in H$,
\begin{eqnarray*}
    DF[a](h, h) &=& \sum_{j = 1}^k \int_{A_j} \left(T[a]^2|\nabla h_j|^2 + 2T[a]DT[a](h) \nabla a_j\nabla h_j\right) dx\\
    &\geq& \sum_{j = 1}^k \left[\int_{A_j} \left(T[a]^2|\nabla h_j|^2\right)dx - \int_{A_j} |2T[a]DT[a](h)\nabla a_j\nabla h_j| dx\right]\\
    &\geq& \lambda^2 \left(\|h\|_{H}^2 - \sum_{j = 1}^k \frac{\Lambda}{\lambda^2} \|DT[a](h)\|_{L^4(A_j)}\|\nabla a_j\|_{L^4(A_j)}\|\nabla h_j\|_{L^2(A_j)}\right).
\end{eqnarray*}
It follows from the continuous embedding of $H^1(D)$ into $L^4(D)$
and (\ref{2.5}) that
\[
    DF[a](h, h) \geq \lambda^2\left(1 - \frac{C_D\Lambda^2
    \theta}{\lambda^2}\right)\|h\|_{\Hd}^2,
\] and therefore, inequality (\ref{5.4.2}) holds true. \qed

We now make use of Theorem \ref{Lemma F derivative} in order to
prove a local Landweber condition which guarantees the convergence
of the reconstruction algorithm.

Let $a$ and $a'$ be in $K$. We can find $t \in [0, 1]$ such that
\begin{equation}
    \|F[a] - F[a']\|_{H^*} = \|DF[ta + (1 - t)a'](a - a')\|_{H^*}
                            \geq C\|a - a'\|_H
\label{Stability}
\end{equation}
by (\ref{5.4.2}). Hence, if $\|a - a'\|_H$ is small enough,  then
\begin{equation}
    \|F[a] - F[a'] - DF[a](a - a')\|_{H^*} \leq \eta \|F[a] - F[a']\|_{H^*}
\label{Landweber condition}
\end{equation}
for some $\eta < \frac{1}{2}.$ In other words, $F$ satisfies the local
Landweber condition (see \cite{Hankeetal:nm1995}).

\subsection{Landweber iteration}

Going back to equation (\ref{5.1}), we have
\begin{equation}
    \Div \Phi^2\nabla a = \Delta \psi
    \label{eqn 2nd}
\end{equation} in the sense of distributions. However, the equation above can be
understood in the classical sense in each inclusion $A_i$. This
observation plays an important role in reconstructing the true
coefficient from the initial guess given in Subsection
\ref{sec:simple}.

Considering $\Delta\psi$ as an element of $H^*$ defined by
\begin{equation}\label{defpsif}
    -\Delta \psi (v) = \sum_{j = 1}^k \int_{A_j} \nabla \psi \cdot \nabla v_j dx,
\end{equation}  for all $v = (v_1, \cdots, v_k),$ we rewrite (\ref{eqn 2nd})
as
\begin{equation}
    F[a] = \Delta \psi.
\label{6.11}
\end{equation}
%
%
%

Recalling that $K$ is closed and convex in $H$, we can employ the classical Hilbert
projection theorem to define the projection from $H$ onto $K$ as
\begin{equation}
    P: H \ni h \mapsto \mbox{argmin}\{\|h - a\|_{H}: a \in K\}.
\end{equation}
It is not hard to verify that
\begin{equation}
    \|P(h) - a\|_H \leq \|h - a\|_H
\label{projection property}
\end{equation} for all $a \in K.$

We next solve (\ref{6.11}) using the Landweber method to minimize
\[
    I(a) = \frac{1}{2} \|F[a] - \Delta \psi\|_{H^*}^2,
\] where $a$ varies in $K$ with the initial guess $a_I = (\alpha_1, \cdots, \alpha_k),$ obtained in Subsection \ref{sec:simple}. The corresponding guess for the coefficient is
\[
    a_I =\sum_{i = 1}^k\alpha_i\1_{A_i}.
\]

There is a gap if we minimize $I$ by the classical Landweber
sequence given by
\[
\begin{array}{rcl}
    a^{(0)} &=& a_I,\\
    a^{(n + 1)} &=& a^{(n)} - \mu DF[a^{(n)}]^*(F[a^{(n)}] - \Delta \psi)
\end{array}
\]
because $a^{(1)}$ may not belong to $K$ and $F[a^{(1)}]$ is not well-defined. Motivated by (\ref{projection property}), which implies $P(a^{(n)})$ is closer to $a_*$ than $a^{(n)}$ is, we modify this formula as
\begin{equation}
    a^{(n + 1)} = P(a^{(n)}) - \mu DF[P(a^{(n)})]^*(F[P(a^{(n)})] - \Delta \psi).
\label{6.4444}
\end{equation}
We have the following convergence result.
\begin{theo} \label{theocv}
Suppose that the true optical distribution $a_* \in K$. Let
$a^{(n)}$ be defined by (\ref{6.4444}) with $a^{(0)}$ being the
initial (piecewise constant) guess obtained as the minimizer of
(\ref{discfunc}). Then the sequence $a^{(n)}$ converges in $H$ to
$a_*$ as $n \rightarrow\infty$.
\end{theo}

Noting that $F$ satisfies the local Landweber condition (see (\ref{Landweber condition})),
we can repeat the proof of Proposition 2.2 in \cite{Hankeetal:nm1995} to see that
\[
    \|a^{(n + 1)} - a_*\|_H^2 + (1 - 2\eta)\|F[P(a^{(n)}) - \Delta \psi]\|^2_{H^*} \leq \|P(a^{(n)}) - a_*\|_H^2.
\]
This and (\ref{projection property}) imply
\begin{equation}
    \|P(a^{(n + 1)}) - a_*\|_H^2 - \|P(a^{(n)}) - a_*\|_H^2 \leq (2\eta - 1)\|F[P(a^{(n)}) - \Delta \psi]\|^2_{H^*} \leq 0.
    \label{6.4}
\end{equation}
It follows that
\[
    \sum_{i = 1}^{\infty}\|F[P(a^{(n)})] - \Delta \psi\|_{H^*}^2 \leq \frac{1}{1 -
    2\eta}\|a_*\|_H^2,
\] and hence
\begin{equation}
    F[P(a^{(n)})] \rightarrow \Delta \psi \mbox{ in } \h \mbox{ as } n \rightarrow \infty.
\label{6.5}
\end{equation}

On the other hand, we can see from (\ref{6.4}) that the sequence
$(P(a^{(n)}))_{n \geq 1}$ is bounded in $H$. Assume that
$P(a^{(n)})$ converges weakly to $a'$ for some $a' \in H$. Since
$K$ is closed and convex, it is weakly closed and therefore $a'
\in K.$ Passing to a subsequence if necessary, this sequence
converges to $a'$ {\it a.e.} and also converges strongly to $a'$
in $\prod_{j = 1}^{k}L^2(A_j).$ So, $T[P(a^{n})]$ converges to
$T[a']$ in $H^1(\Omega)$ and hence in $L^4(\Omega).$ For all $v
\in H,$ we have $$\begin{array}{l}
    \ds \sum_{j = 1}^k \int_{A_j}(T[P(a^{(n)})]^2\nabla P(a^{(n)}) - T[a']^2\nabla a')\nabla v
    dx\\
    \nm \ds
   = \sum_{j = 1}^k \bigg[\int_{A_j}(T[P(a^{(n)})]^2 - T[a']^2)\nabla P(a^{(n)})\nabla v dx + \int_{A_j}T[a']^2(\nabla P(a^{(n)}) - \nabla a')\nabla v dx\bigg],
\end{array}$$
which goes to $0$ by the dominated convergence theorem and the weak convergence
of $P(a^{(n)})$ to $a'$ in $H$. We have obtained
$
    F[a'] = \Delta \psi = F[a_*].
$ Using (\ref{5.4.2}) gives $a' = a_*.$

In summary, if the true coefficient $a_*$ is a perturbation of a
constant on each inclusion then the coefficient $a_I$ obtained in
Section \ref{sec:simple} is quite closed to $a_*$. Moreover, the
misfit between the initial guess $a_I$ and the true distribution
$a_*$ can be properly corrected by the sequence in (\ref{6.4444}).

\section{Concluding remarks} In this paper we have introduced a Landweber
scheme for reconstructing piecewise smooth optical absorption
distributions from opto-acoustic measurements and proved its
convergence. Because of the jumps in the absorption coefficient,
we have used weak formulations for the Helmholtz decomposition for
$\Phi_*^2\nabla a_*$ and the relation between  the spherical Radon
transform of its gradient part $\psi$ and the cross-correlation
measurements $M_\eta(y, r)$. Note that we can enrich the set of
data as follows. For $f \in L^2(\partial \Omega)$ such that $f
\geq 0$ {\it a.e.} on $\partial \Omega$, compute instead of
(\ref{cross}) the quantity
$$
  M^{f,g}_\eta(y, r) = \frac{1}{\eta^2}\int_{\partial \Omega} (f \partial_\nu
   \Phi^g_{u_{y,r}^\eta} - g \partial_\nu \Phi^f)\, d\sigma, \hspace*{.24in}
    y \in S_{\mu}, r > 0.
$$
Similarly to (\ref{data}), integration by parts yields
\begin{equation} \label{dataenr}
    M^{f,g}_\eta(y, r) = \frac{1}{\eta^2}\int_\Om (a_{u_{y, r}^\eta} - a)\Phi^f \Phi^g_{u_{y, r}^\eta}dx,
\end{equation}
where $\Phi^f$ is the solution of (\ref{Optical eqn}) with $g$
replaced by $f$.

The enriched data (\ref{dataenr}) may be used in order to
generalize our approach to the case of
 measurements of the outgoing light intensities on only part of
$\partial \Omega$ by choosing $f$ supported only on the accessible
part of the boundary. Another interesting and challenging problem
is to prove statistical stability of the proposed reconstruction
with respect to a measurement noise by combining Fourier
techniques together with statistical tools \cite{ip12}. Numerical
implementation of the Landweber-type iteration is under
consideration and will be the subject of a forthcoming
publication. The behavior of the proposed method with respect to
the optical absorption contrast will be investigated. It is
expected that higher the contrast, more efficient the method is.

\appendix

\section{Proof of Lemma \ref{The unique continuation property}}
\label{appenA}

    The boundedness of $\phi$ together with the assumption that $\phi \equiv 0$ on $\partial \Omega$ imply  by standard
    regularity results
    that $\phi \in
    \C^1( \partial \Omega \cup \Omega \setminus \overline \Omega')$.
    Arguing similarly to Proposition 2.5 in \cite{AmmariBossyLocLaurentetal:preprint2012},
    we see that $\phi \in \C^2(\Omega \setminus \overline \Om').$ Define
    \[
        \mathcal{U} = \{x \in \Omega \setminus \overline \Om': u(x) \not = 0\}.
    \] The continuity of $\phi$ shows that $\mathcal{U}$ is open. Assume, on contrary, that $\mathcal{U}$ is nonempty.

    Noting that $\mathcal{U}$ can be decomposed as the union of its connected open subsets. Denote by $\mathcal{O}$
     the connected component of $\mathcal{U}$, which is closest to $\partial \Omega.$ Without loss
     of generality, assume that $\phi > 0$ in $\mathcal{O}.$ Let
    \[
        \delta = \mbox{dist}(\mathcal{O}, \partial \Omega).
    \] The distance above is understood as the length of the shortest curve, contained
    in $\Omega \setminus \overline \Om'$ and connecting $\overline {\mathcal{O}}$ and $\partial \Omega.$

    In the case that $\delta = 0,$ $\partial \mathcal{O}$ and $\partial \Omega$ have a
    common point $x_0$. Applying the Hopf lemma for the equation
    \[
        \left\{
            \begin{array}{rcll}
                -\Delta \phi + c\phi &=& 0 &\mbox{in } \mathcal{O},\\
                \nm
                \phi &>& 0 &\mbox{on } \partial \mathcal{O},
            \end{array}
        \right.
    \]
gives $\partial_{\nu}\phi(x_0) < 0,$ which is impossible.

When $\delta > 0,$ it is easy to see that $\phi \equiv 0$ in a
neighbourhood of $\partial \Omega$. Assume that such a neighbourhood
and $\mathcal{O}$ have a common boundary point $x_0$. Noting that
$\nabla \phi(x_0) = 0$, we can apply the Hopf lemma again to get
the contradiction. \cqfd

\section{Proof of Proposition \ref{approx}}
\label{appenB} We write $u$ and $v$ when referring to
$u_{y,r}^\eta$ and $v_{y,r}^\eta$ respectively for simplicity.
Using (\ref{2.4'}) and the same arguments when estimating $\|a_u -
a\|_{L^1(\Omega)}$ in the previous section yields
\[
    \|a_{u} - a\|_{L^2(\Omega)} \leq O(\eta).
\] This, together with standard $H^2$-regularity results (see, for instance,  \cite[Theorems 8.8 and
8.12]{GilbargTrudinger:1977}) and the embedding of $H^2(\Omega)$
into $L^{\infty}(\Om)$, gives
\[
    \|\Phi_u - \Phi\|_{L^{\ii}(\Omega)} \leq O(\eta).
\]
Hence, it follows from (\ref{3.5}) that
\begin{eqnarray}
\left|\int_\Om(a_u-a)\Phi\Phi_u dx - \int_\Om(a_u - a)\Phi^2
dx\right| &\leq&
\int_\Om\Phi|a_u-a||\Phi_u-\Phi|dx \nonumber \\
&\leq&
\n\Phi{L^\ii(\Omega)}\n{a_u-a}{L^1(\Omega)}\n{\Phi_u-\Phi}{L^\ii(\Omega)}
\nonumber
\\ &\leq& c\eta^3. \label{3.3}
\end{eqnarray}
The constant $c$ depends only on $ \overline{a} = \max a$ and
$\underline{a} = \min a$, both of which are assumed to be known.
The independence of $c$ on $||\Phi||_{L^\infty(\Omega)}$ can be
deduced from Lemma \ref{lemma: varphi lambda Lambda}. Now, note
that the second integral in the left hand side of (\ref{3.3}) can
be rewritten as
\begin{equation}\nonumber\begin{aligned}
\int_\Om(a_u-a) \Phi^2 dx &=\sum_{i=1}^n\int_{A_i\cup P(A_i)}(a_u-a)\Phi^2 dx\\
&=\sum_{i=1}^n\left[\int_{A_i\cap P(A_i)}(a_u-a)\Phi^2 dx
+\int_{A_i\La P(A_i)}(a_u-a)\Phi^2 dx\right],
\end{aligned}\end{equation}
and that the integral in (\ref{tilde M}) is equal to
\begin{equation}\nonumber\begin{aligned}
\int_\Om (a-a_0)\Div(\Phi^2 v) dx &= \sum_{i=1}^n\int_{A_i} (a-a_0)\Div(\Phi^2 v)dx\\
&=\sum_{i=1}^n\left[\int_{\dr\A_i}(a_i-a_0)\Phi^2v \cdot \nu_i d\sigma - \int_{A_i} \Phi^2\g a \cdot v dx\right]\\
&= \sum_{i=1}^n\left[\int_{\dr\A_i}(a_i-a_0)\Phi^2v \cdot \nu_i
d\sigma \right. \\ & \hspace*{.60in} - \left. \int_{A_i\cap
P(A_i)} \Phi^2\g a\cdot v dx - \int_{A_i\bs P(A_i)} \Phi^2\g
a\cdot vdx\right].
\end{aligned}\end{equation}
Therefore, we have
\begin{eqnarray*}
&& \left|\int_\Om(a_u-a)\Phi^2dx + \int_\Om (a-a_0)\Div(\Phi^2 v)dx\right| \\
&& \hspace*{.24in} \leq \sum_{i=1}^n \left|\int_{A_i\cap
P(A_i)}(a_u-a+ \g a \cdot v)\Phi^2dx\right|
 + \left|\int_{A_i\bs P(A_i)}\Phi^2\g a_i \cdot v dx\right|\\
&& \hspace*{.48in} + \left|\int_{A_i\La P(A_i)}(a_u-a)\Phi^2 dx
-\int_{\dr\A_i}(a_i-a_0)\Phi^2v \cdot \nu_i d\sigma\right|.
\end{eqnarray*}
Denote by  $\al_i$, $\be_i$ and $\gamma_i$ the last three
quantities in the inequality above. We need to prove that they all
are bounded by $O(\eta^3)$ to complete the proof.
\begin{itemize}
\item[(i)] Since for all $i = 1, \cdots, k$, $a_i \in
\C^2(\overline A_i)$ and $\|D^2 a_i\|_{L^{\infty}(A_i)}$ are
bounded by some known constants, we can find a constant $c$  such
that
\begin{equation}\nonumber\begin{aligned}
\left|\int_{A_i\cap P(A_i)}(a_u-a-\g a \cdot u)\Phi^2dx\right|
\leq \n{\Phi}{L^\ii(A_i)}^2\n{D^2a}{L^\ii(A_i)}\eta^2|\Sigma_\eta|
\leq c^i_1\eta^3.
\end{aligned}\end{equation}
On the other hand, using the classical substitution method in
integration gives
\begin{eqnarray*}
&& \hspace*{-.24in} \left|\int_{A_i\cap P(A_i)\cap\Sigma_\eta}\g a \cdot (u+v)\Phi^2 dx\right|\\
&&\leq\left|\int_{S(y,r)}(\Phi^2\g a)(\xi) \cdot \int_{-\eta}^\eta\left(\1_{A_i}(u+v)\right)\left((1+\frac\rho r)\xi\right)d\rho d\xi\right| \\
&& \hspace*{.24in} + \left|\int_{S(y,r)}\int_{-\eta}^\eta\left[(\Phi^2\g a)\left((1+\frac\rho r)\xi\right)-(\Phi^2\g a)(\xi)\right].\left(\1_{A_i}(u+v)\right)\left((1+\frac\rho r)\xi\right)d\rho d\xi\right|\\
&& \leq 0 + 2\eta^2\n{\dr_r(\Phi^2\g a)}{L^\ii(A_i)}|\Sigma_\eta| \\
&& \leq c \eta^3.
\end{eqnarray*}
The quantity $\alpha_i$ is bounded from above by $O(\eta^3)$
because
\begin{equation}\nonumber\begin{aligned}
\al_i &\leq \left|\int_{A_i\cap P(A_i)}(a_u-a-\g a \cdot u)\Phi^2 dx \right|+\left|\int_{A_i\cap P(A_i)\cap\Sigma_\eta}\g a \cdot (u+v)\Phi^2 dx\right|.\\
\end{aligned}\end{equation}

\item[(ii)] The fact that $\beta_i \leq O(\eta^3)$ can be deduced
from the boundedness of the integrand and (\ref{2.4'}).

\item[(iii)] The main point of the proof is the estimate of
$\gamma_i$. Denote $\ti v(z)={v(z)}/{|v(z)|}$ when it is defined
and
\begin{equation}\nonumber\begin{aligned}
\chi:\dr\A_i\times[0,\eta]\fg A_i\La P(A_i)\\
(z,t)\longmapsto z+t\ti v(z).
\end{aligned}\end{equation}
$\chi$ is well defined when $v(z)$ is non-zero and not parallel to
$\dr A_i$. For any $z\in \dr\A_i$ satisfying this condition, we
denote $T(z)$ the tangent plane to $\dr A_i$ and $B(z)$ a basis
adapted to the sum $\R^d=T(z)\oplus\R\ti v(z)$. Let $I_{d-1}$
denote the $(d-1) \times (d-1)$ identity matrix. We get
$$ d\chi(z,t)=
 \begin{bmatrix}
  I_{d-1}+td\ti v(z)  & 0 \\
  *  & 1
 \end{bmatrix}
$$
and as $\ti v(z)=({z-y})/{|z-y|}$. The operator $d\ti v(z)$ does
not depend
 on $\eta$ and  $td\ti v(z)= O(\eta)$ with a constant depending on $r_0$. Then,
$$\det(d\chi(z,t))=1+t\ \Div(\ti v)(z)+ O(\eta^2)=1+ O(\eta).$$
As $B(z)$ is not orthonormal, the differential volume written with
the variables $(z,t)$ depends on the angle between $\ti v(z)$ and
$\nu(z)$ called $\theta(z)$. This volume at the point $z+t\ti
v(z)$ is $(1+ O(\eta))\cos(\theta(z))dtdz$. Knowing this, we
denote
$$(\dr A_i)^\pm=\{z\in\dr A_i,\ \pm\theta(z)>0\}$$
and write
$$
\begin{array}{lll}
\ds \int_{P(A_i)\bs A_i}(a_u-a)\Phi^2 dx &=&\ds \int_{(\dr
A_i)^+}\int_0^{|v(z)|}(a_u-a_0)\Phi^2(z+t\ti v(z))(1+ O(\eta))\\
\nm && \times \cos(\theta(z))dtdz \end{array} $$ and as $a_i$ and
$\Phi$ are $\C^1(\overline A_i)$, we can write that for any
$z\in(\dr A_i)^+$, and $t\in[0,|v(z)|]$,
$$|(a_u-a_0)\Phi^2(z+t\ti v(z))-(a_i-a_0)\Phi^2(z)|\leq O(\eta).$$
Then,
$$\left|\int_0^{|v(z)|}(a_u-a_0)\Phi^2(z+t\ti v(z))(1+ O(\eta))dt-(a_i-a_0)\Phi^2(z)|v(z)|\right|\leq O(\eta^2).$$
\end{itemize}
Now, noticing that $\cos(\theta)|v(z)|=v(z) \cdot \nu(z)$ and that
$\sigma((\dr A_i)^+\cap\Sigma_\eta)$, the surface measure of $(\dr
A_i)^+\cap\Sigma_\eta$, is of order $O(\eta)$, we have
$$\left|\int_{P(A_i)\bs A_i}(a_u-a)\Phi^2-\int_{(\dr A_i)^+}(a_i-a_0)\Phi^2v \cdot \nu\right|\leq O(\eta^3).$$
We also get
$$\left|\int_{A_i\bs P(A_i)}(a_u-a)\Phi^2-\int_{(\dr A_i)^-}(a_i-a_0)\Phi^2v \cdot \nu\right|\leq O(\eta^3)$$ by
the same arguments.\qed

\section{Construction of $\Ra[\psi]$ from  formula (\ref{formulainvert})}
\label{A primitive construction in the sense of distributions}

In order to construct $\Ra[\psi]$ from  formula
(\ref{formulainvert}), we need to invert the operator
$\frac{\dr}{\dr r}:L^2(C)\fg G^{-1}(C)$ and prove the stability of
the inversion.
 For any $f\in L^2(C)$, by Fubini's theorem, the function $F(y,r)=\int_0^rf(y,\rho)d\rho$
 is well-defined and
 in $G(C)$ but not in $G_0(C)$. Since this operator is acting on distributions which are zero on
$S_\mu\times]0,r_0[$, we introduce
\begin{equation}\begin{aligned}\nonumber
p:L^2(C) &\fg G_0(C)\\
\ph &\longmapsto\left[(y,r)\mapsto
-\int_0^r\left(\ph(y,\rho)-\frac R
{r_0}\chi_{]0,r_0[}(\rho)\ph(y,\rho R/r_0)\right)d\rho\right]
\end{aligned}\end{equation}
and its dual
\begin{equation}
p^*:G^{-1}(C) \fg L^2(C).
\end{equation}
The following result holds.
\begin{prop}
For all $f\in L^2(C)$ such that $f=0$ on $S_\mu\times]0,r_0[$, we
have the inversion formula
$$p^* [\frac{\dr f}{\dr r}]=f.$$
\end{prop}
\Proof For any $\ph\in L^2(C)$, we have $\frac{\dr}{\dr
r}p[\ph]=-\ph$ on $S_\mu\times[r_0,R[$ and therefore,
$$\int_Cp^* \left[\frac{\dr f}{\dr r}\right]\ph=
\left<\frac{\dr f}{\dr
r},p[\ph]\right>_{G^{-1}(C),G^1_0(C)}=-\int_Cf\frac{\dr}{\dr
r}p[\ph]=\int_Cf\ph,$$ which yields the claimed result. \qed

\begin{prop}
For all $u\in \M :=\big\{ v \in G^{-1}(C): \mbox{\rm supp}(v)
\subset S_{\mu} \times [r_0, R[ \big\}$,
$$\n {p^*u} {L^2(C)}\leq\n u  {G^{-1}(C)}$$
\end{prop}
\Proof We first note that $p^*[u]=0$ on $S_\mu\times]0,r_0[$.
Then, for any $\ph\in L^2(C)$, we get
\begin{equation}\begin{aligned}\nonumber
\left|\int_Cp^* \left[u\right]\ph\right|=\left|\int_Cp^*
\left[u\right]\chi_{[r_0,R[}\ph\right|
&\leq\n u  {G^{-1}(C)}\n {p[\chi_{[r_0,R[}\ph ]}{G^1_0(C)}\\
&\leq\n u  {G^{-1}(C)}\n {\frac{\dr}{\dr r}p[\chi_{[r_0,R[}\ph ]}{L^2(C)}\\
&\leq \n u  {G^{-1}(C)}\n {\chi_{[r_0,R[}\ph}{L^2(C)}\\
&\leq \n u  {G^{-1}(C)}\n {\ph}{L^2(C)},
\end{aligned}\end{equation} and the proof is complete. \qed

Finally, we deduce the following result.
\begin{cor} From formula (\ref{formulainvert}), we have
$$\Ra[\psi]=\frac{1}{r_0\n w 1} p^*(r^{d-2}M).$$
Moreover, for $\eta$ small, if $\Ra[\psi_\eta]=\frac{1}{r_0\n w 1}
p^*(r^{d-2}M_\eta)$,  then
$$\n {\Ra[\psi-\psi_\eta]}{L^2(C)} \leq\frac{R^{d-2}}{r_0\n w 1}\n{M-M_\eta}{G^{-1}(C)},$$ which insures
the stability of the construction of $\Ra[\psi]$ from the
measurements $M_\eta$.
\end{cor}




\begin{thebibliography}{1}


\bibitem{sergio} {\sc G. Alessandrini and S. Vessella}, {\em Lipschitz stability for
the inverse conductivity problem},  Adv. Appl. Math.,  35 (2005),
pp.~207--241.

\bibitem{AMMARI-08}
{\sc H.~Ammari}, {\em An Introduction to Mathematics of Emerging
Biomedical Imaging}, Vol. {62}, { Mathematics and Applications},
Springer-Verlag, Berlin, 2008.


\bibitem{AMMARI-BONNETIER-CAPDEBOSCQ-08}
{\sc H.~Ammari, E.~Bonnetier, Y.~Capdeboscq, M.~Tanter, and
M.~Fink},
\newblock {\em Electrical impedance tomography by elastic
deformation},
\newblock {SIAM J. Appl. Math.}, 68 (2008), pp.~1557--1573.


\bibitem{AmmariLaurentetal:preprint2012}
{\sc H.~Ammari{, E. Bossy, J. Garnier and L. Seppecher}}, {\em
  Acousto-electromagnetic tomography}, SIAM J. Appl. Math., to
  appear.

\bibitem{AmmariBossyLocLaurentetal:preprint2012}
{\sc H.~Ammari, E. Bossy, J. Garnier, L. H. Nguyen and L.
Seppecher}, {\em A
  reconstruction algorithm for ultrasound-modulated diffuse optical
  tomography}, Proc. Amer. Math. Soc., submitted.



\bibitem{bretin} {\sc H. Ammari, E. Bretin, V. Jugnon, and A. Wahab},
{\em Photoacoustic imaging for attenuating acoustic media},
Lecture Notes in Math., Vol. 2035, pp.~57--84, Springer-Verlag,
Berlin, 2011.

\bibitem{siap2011} {\sc H.~Ammari, Y.~Capdeboscq, F. de Gournay, A. Rozanova, and F.
Triki}, {\em Microwave imaging by elastic perturbation}, {SIAM J.
Appl. Math.}, 71 (2011), pp.~2112--2130.


\bibitem{ip12} {\sc H.~Ammari, J. Garnier, and W. Jing}, {\em Resolution and stability analysis in acousto-electric
imaging}, Inverse Problems, to appear.

\bibitem{simon} {\sc S. R. Arridge}, {\em Optical tomography in medical
imaging}, Inverse Problems, 15 (1999),  R41--R93.

\bibitem{beretta} {\sc E. Beretta and E. Francini},
{\em Lipschitz stability for the electrical impedance tomography
problem: the complex case}, {Comm. Partial Differential
Equations},  {36} (2011), pp.~{1723--1749}.

\bibitem{beretta2} {\sc E. Beretta, M. V. De Hoop, and L. Qiu},
{\em Lipschitz stability of an inverse boundary value problem for
a Schr\"odinger type equation}, preprint.

\bibitem{born} {\sc M. Born and E. Wolf}, {\em Principles of Optics},  Cambridge University
Press, Cambridge, 1999.

\bibitem{CAPDEBOSCQ-FEHRENBACH-DEGOURNAY-KAVIAN-09}
{\sc Y.~Capdeboscq, J.~Fehrenbach, F.~de~Gournay, and O.~Kavian},
\newblock {\em Imaging by modification: numerical reconstruction of local
  conductivities from corresponding power density measurements},
\newblock {SIAM J. Imag. Sci.}, 2 (2009), pp.~1003--1030.

\bibitem{galdi} {\sc G.~P. Galdi}, {\em An Introduction to the
Mathematical Theory of the Navier-Stokes Equations, Vol. I,
Linearized Steady Problems}, Springer-Verlag, New York, 1994.

\bibitem{GEBAUER-SCHERZER-08}
{\sc B.~Gebauer and O.~Scherzer}, {\em Impedance-acoustic
tomography}, {SIAM J. Appl. Math.}, 69 (2008), pp.~565--576.

\bibitem{GilbargTrudinger:1977}
{\sc D.~Gilbarg and N.~S. Trudinger}, {\em Elliptic partial differential
  equations of second order}, Springer-Verlag, Berlin Heidelberg New York,
  1977.

\bibitem{Hankeetal:nm1995}
{\sc M.~Hanke, A. Neubauer, and O. Scherzer}, {\em A convergence
analysis of
  the Landweber iteration for nonlinear ill-posed problems}, Numer. Math., 72
  (1995), pp.~21--37.

\bibitem{otmar} {\sc M. V. de Hoop , L. Qiu, and O. Scherzer}, {\em Local
analysis of inverse problems: H\"older stability and iterative
reconstruction}, Inverse Problems, 28 (2012), 045001.



\bibitem{john2} {\sc K. Kilgore, S. Moskow, and J.~C. Schotland}, {\em Inverse
Born series for diffuse waves},  Imaging microstructures,
pp.~113--122, Contemp. Math., 494, Amer. Math. Soc., Providence,
RI, 2009.

\bibitem{LadyzhenskayaUraltseva:sv1985}
{\sc O.~A. Ladyzhenskaya and {N.N. Ural'tseva}}, {\em Linear and Quasilinear
  Elliptic Equations}, Translated from the Russian by Scripta Technica, Inc.
  Translation editor: Leon Ehrenpreis, Academic Press, New York - London, 1968.


\bibitem{landweber} {\sc L. Landweber}, {\em An iteration formula for Fredholm integral equations of the
first kind}, American J. Math., 73 (1951),  pp.~615--624.

\bibitem{LocSchmitt:die2009}
{\sc N.~H. Loc and K.~Schmitt}, {\em On positive solutions of quasilinear
  elliptic equations}, Differential and Integral Equations, 22 (2009),
  pp.~829--842.



\bibitem{born3} {\sc V. A. Markel and J. C. Schotland}, {\em Symmetries, inversion formulas,
and image reconstruction for optical tomography}, Phys. Rev. E, 70
(2004), 056616.

\bibitem{john3} {\sc S. Moskow and J.~C. Schotland}, {\em
Convergence and stability of the inverse scattering series for
diffuse waves}, Inverse Problems, 24 (2008), 065005.


\bibitem{rossum} {\sc M.~C.~W. van Rossum and Th.~M.
Nieuwenhuizen}, {\em Multiple scattering of classical waves:
microscopy, mesoscopy, and diffusion}, Rev. Modern Phys., 71
(1999), pp.~313--371.


\bibitem{papa} {\sc L. Ryzhik, G. Papanicolaou, and J.~B. Keller},
{\em  Transport equations for elastic and other waves in random
media}, Wave Motion 24 (1996),  pp.~327--370.


\bibitem{john} {\sc J. C. Schotland}, {\em Direct reconstruction
methods in optical tomography}, Lecture Notes in Math., Vol. 2035,
pp.~1--29, Springer-Verlag, Berlin, 2011.


\bibitem{john4} {\sc J. C. Schotland and V.~A. Markel},
{\em Inverse scattering with diffusing waves}, J. Opt. Soc. Amer.
A, 18 (2001), pp.~2767--2777.


\bibitem{otmar2} {\sc T. Widlak and O.
Scherzer}, {\em Hybrid tomography for conductivity imaging},
Inverse Problems, to appear.

\end{thebibliography}
\end{document}